\documentclass[APA,LATO1COL]{WileyNJDv2}

\articletype{Article Type}%

\received{26 April 2016}
\revised{6 June 2016}
\accepted{6 June 2016}

\usepackage{amsfonts,amsmath,amsthm,amssymb,graphicx,mathrsfs}

\newtheorem{lem}{Lemma}[section]
\newtheorem{cor}{Corollary}[section]

\newtheorem{theo}{Theorem}[section]

\begin{document}

\title{On a low-rank matrix single index model}

\author{The Tien Mai}

\authormark{The Tien Mai}

\address{
	\orgdiv{Department of Mathematical Sciences},
	\orgname{Norwegian University of Science and Technology}, 
	\orgaddress{\state{Trondheim}, \country{Norway}}}

\corres{\email{the.t.mai@ntnu.no}}


\abstract[Abstract]{
In this paper, we present a theoretical study of a low-rank matrix single index model. This model is recently introduced in biostatistics however its theoretical properties on estimating together the link function and the coefficient matrix are not yet carried out. Here, we advance on using PAC-Bayesian bounds technique to provide a rigorous theoretical understanding for jointly estimation of the link function and the coefficient matrix.   
}


\maketitle


\section{Introduction}

In this work, we study the following single index model where the response variable $ Y\in \mathbb{R} $ and the covariate matrix $ X \in \mathbb{R}^{d\times d} $ satisfy
\begin{equation}
\label{main_model}
Y = f^*( \left\langle X,B^* \right\rangle ) + \epsilon
\end{equation}
where $ \left\langle X,B \right\rangle = {\rm trace}(XB) $ and $ B \in \mathbb{R}^{d\times d} $ is an unknown coefficient matrix. The link function $ f^* $ is an unknown univariate measurable function. The noise term $ \epsilon $ is assumed to have $ 0 $-mean. Here, $ B $ is a low-rank and symmetric matrix. For model identifiability, we further assume that the Frobenius norm of $ B $ satisfies that $ \| B \|_F =1 $. 

A similar model, where the unknown matrix $ B $ is assumed to be sparse, is recently used by \cite{weaver2021single} in biostatistics to study the relationship between a response variable and the functional connectivity associated with a brain region. More recently, this model has been studied in \cite{fan2022understanding} where the authors focus on estimating the unknown low-rank matrix $ B $ through implicit regularization. This model can be seen as a nonparametric version of the trace regression in \cite{koltchinskii2011nuclear,zhao2017trace}. The trace regression, in which the link function is the identity, covers a wide range of statistical models such as reduced rank regression, matrix completion and linear regression.

The single index model is a particularly useful generalization of the linear model \cite{nelder1972generalized} as it can be interpreted naturally that the model changes only in the direction of the (vector/matrix) parameter and this change is described by the link function $ f^\star $. Various works on applications and extension had been done for vector single index model, see for example \cite{kong2007variable,alquier2013sparse,weaver2021single,putra2020study}. 

Let $ S_{1,+}^d $ denote the set of all symmetric matrix $ B \in \mathbb{R}^{d\times d} $ such that $ \| B \|_F =1 $. We define the expected risk, for any measurable $ f:\mathbb{R}\rightarrow \mathbb{R}  $ and $ B \in S_{1,+}^d $, as
\begin{equation*}
R(B,f) = \mathbb{E} \left[ (Y- f(\left\langle X,B \right\rangle )^2) \right]
\end{equation*}
and denote the empirical counterpart of $ R(f,B) $ by
\begin{equation*}
r_n( B, f ) = \frac{1}{n} \sum_{i=1}^{n}
(Y_i- f(\left\langle X_i,B \right\rangle )^2).
\end{equation*}
In this work, we study the predictive aspects of the model: a pair $ (f,B) $ predicts almost as well as $ (f^\star,B^\star) $ if $ R( B,f)- R ( B^\star,f^\star) $ is small.

Our approach is based on PAC-Bayesian bound technique \cite{catoni2007pac} which is a useful tool to obtain oracle inequalities bound. It is also an important point, similar to Bayesian analysis, in PAC-Bayesian bound is to specify prior distribution over the parameter space. Our prior distribution for the link function is from the reference \cite{alquier2013sparse}, while the prior distribution is inspired from the eigen decomposition of the matrix parameter $ B $. The details of our approach and the prior specifications are given in the next section.

\section{Main result}

\subsection{Method}
We assume that, in model \eqref{main_model}, $\mathbb E[ \epsilon | X] = 0 $ and, consequently, that $\mathbb E \epsilon =0 $. We will
put the following assumption on the distribution of $\epsilon $. 

\noindent\textbf{Assumption $\mathbf N$.}  There exist two positive constants
$\sigma$ and $L$ such that, for
all integers $k \geq 2$, 
$$ 
\mathbb{E}\left[|\epsilon|^{k}\,|\, X\right] 
\leq
\frac{k!}{2}\sigma^{2}L^{k-2}.
$$

\noindent The above assumption means that the noise follows a sub-exponential distribution which includes for example Gaussian noise or bounded noise, see \cite{boucheron2013concentration}.

We will also need that the random variable $X $ is
almost surely bounded by a constant. Moreover, the link function $f^{\star} $ is also assumed to be bounded by some known positive constant. Thus, denoting by
$\|X \|_{\infty}$ the supremum norm of $X $ and by $\|f^{\star}\|_{\infty}$
the functional supremum norm of $f^{\star}$ over $[-1,1]$, we set the following assumption.

\noindent\textbf{Assumption $\mathbf B$.} The condition $\|X \|_{\infty} \leq 1$ holds almost surely and there exists a positive constant $C$ larger than 1 such that $\|f^{\star}\|_{\infty} \leq C.$

\noindent To make the technical proofs as clear as possible, no attempt was done to optimize the constants. In particular, the requirement $C\geq1$ is completely technical and it could be removed by replacing $C$ by $\max[C,1]$ throughout the paper.

In order to approximate the link function $f^{\star}$, we consider the vector
space $\mathcal F$ spanned by a given countable dictionary of measurable
functions $\{\varphi_{j}\}_{j=1}^{\infty}$. Following \cite{alquier2013sparse}, the approximation
space $\mathcal F$ is the set of (finite) linear combinations of functions of
the dictionary. Each $\varphi_j$ of the collection is assumed
to be defined on $[-1,1]$ and to take values in $[-1,1]$. For simplicity, we assume that each $\varphi_j$ is differentiable and such that, for some positive constant $ C_\phi $, 
$$
\|\varphi_j'\|_\infty \leq j C_\phi   .
$$
For example, the (non-normalized) trigonometric system
$\varphi_1(t)=1,\, \varphi_{2j}(t)=\cos(\pi j t),\,
\varphi_{2j+1}(t)= \sin (\pi j t), j=1, 2, \hdots $ satisfies this assumption.

Our approach is inspired from the work by \cite{alquier2013sparse} where the authors explored the PAC-Bayesian approach in \cite{catoni2007pac} for sparse-vector single index model. The method need first to specify a probability measure $\pi$ on $\mathcal S_{1,+}^d \times \mathcal F$, called the prior, which in our framework should enforce the  properties of the target regression function and parameter. In this work, we first let
$$ \mbox{d}\pi( B,f) = \mbox{d}\mu( B) \mbox{d}\nu(f),$$
in other words, it means that the distribution over the index matrix is assumed to be independent of the distribution over the link functions.

As $ B $ is a symmetric matrix, it can be represented in eigen-decomposition form as
$
B = U \Lambda U^\top
$
where $ U $ is a orthogonal matrix $ UU^\top = UU^{-1} = I_d $ (identity matrix of dimension $ d\times d $ ) and $ \Lambda $ is the diagonal matrix whose diagonal elements are the corresponding eigenvalues $ \lambda_1, \ldots, \lambda_d $. To assure that $ \| B \|_F =1 $, we need to have that $ \sum_{j=1}^{d}\lambda_i^{2} =1 $. This is because we have $ \| B \|_F = \sqrt{{\rm trace}(B^2)} $ and $ {\rm trace}(B^2) = \sum_{j=1}^{d}\lambda_i^{2} $. Moreover, we are also interested in imposing low-rankness on B. Thus it is naturally assuming that most of $ \lambda_1, \ldots, \lambda_d $ are quite close to zeros and just a few of them are significantly larger than the rest.

\noindent With this observation in mind, the prior for $ B $ is proposed as follows. We simulate an orthogonal matrix $ V $ and simulate $ (\gamma_1,\ldots, \gamma_d) $ from a Dirichlet distribution $ Dir(\alpha_1,\ldots,\alpha_d) $. Put 
$$ 
B = V {\rm diag} (\gamma_1^{1/2},\ldots, \gamma_d^{1/2})V^\top .
$$
We assume that $ \sum_{i=1}^{d} \alpha_i =1 $ and that $ \alpha_i \geq D_1 $.
In order to get an approximate low-rank matrix, one can take all parameters of the Dirichlet distribution equal to a constant that is very closed to $ 0 $, for example we focus on $ \alpha_1 = \ldots = \alpha_d = 1/d $. It is noted that a typical drawing will lead to one of the $ \gamma_i $s close to 1 and the others close to 0. See \cite{wallach2009rethinking} for discussions on choosing the parameters for Dirichlet distribution.

Now, we present a prior distribution on $ \mathcal F $ that we have opted to use the prior introduced in \cite{alquier2013sparse}. Let us define, for any positive integer $M\leq n$ and all $\Lambda>0$, 
\begin{equation*}
\mathcal B_{M}(\Lambda) 
=
\left\{(
\beta_{1},\hdots,\beta_{M})\in\mathbb{R}^{M}:
\sum_{j=1}^M j|\beta_{j}| \leq \Lambda \text{ and } \beta_{M} \neq 0 \right\}.
\end{equation*}
Now, let $\mathcal F_M(\Lambda)\subset \mathcal F$
denote the image of $\mathcal B_M(\Lambda)$ by the map 
\begin{equation*}
\begin{array}{clcl}
{\Phi}_{M}& :\mathbb R^M &\to & \mathcal F\\
& (\beta_{1},\hdots,\beta_{M}) &\mapsto &\sum_{j=1}^{M} \beta_{j}\varphi_{j}.
\end{array}
\end{equation*}
It is noted that Sobolev spaces are well approximated by $\mathcal F_M(\Lambda)$  as $M$ grows. Finally, we define $\nu_{M}(\mbox{d}f)$ on the set $\mathcal
F_M(C+1)$ as the image of the
uniform measure on $\mathcal B_M(C+1)$ induced by the map ${\Phi}_{M}$. The choice of the prior $\nu$ on $\mathcal F$ is as
\begin{equation}
\label{2plus}
\mbox{d} \nu(f) 
=  
\frac{\displaystyle \sum_{M=1}^{n}10^{-M}
	\mbox{d}\nu_{M}(f)}{1-(\frac{1}{10})^{n}} .
\end{equation}

The reason of choosing $C+1$ instead of $C$ in the definition of the prior support is essentially technical. This bound ensures that when the target $f^{\star}$ belongs to $\mathcal F_n(C)$, then a small ball around it is contained in $\mathcal F_n(C+1)$. It could be safely replaced by $C+u_n$, where $\{u_n\}_{n=1}^{\infty}$ is any positive sequence vanishing sufficiently slowly as $n\to \infty$. The integer $M$ can be viewed as a measure of the
``dimension'' of the function $f$---the larger $M$, the more complex the
function---and the prior $\nu$ adapts again to the sparsity idea
by penalizing large-dimensional functions $f$. The coefficient $10^{-M}$ which appears in \eqref{2plus} shows that more complex models have a geometrically decreasing influence. Inspiring from practical results in \cite{alquier2013sparse}, the value 10 is somehow arbitrary. It could be in general replaced by a coefficient $\alpha$ at the price of a more technical analysis. 

Now, let $\lambda$ be a positive real number, called the inverse temperature parameter
hereafter. The estimates $\hat{ B}_{\lambda}$ and $\hat{f}_{\lambda}$ of
$ B^{\star}$ and $f^{\star}$, respectively, are simply obtained by randomly drawing
$$ (\hat{ B}_{\lambda},\hat{f}_{\lambda}) \sim \hat{\rho}_{\lambda},$$
where $\hat{\rho}_{\lambda}$ is the so-called Gibbs posterior distribution over $\mathcal S_{1,+}^p\times \mathcal F_n(C+1)$, defined by the probability density 
$$ \frac{\mbox{d}\hat{\rho}_{\lambda}}{\mbox{d}\pi}( B, f)= \frac{\exp\left
	[-\lambda r_n( B,f)\right]}{\displaystyle \int \exp\left [-\lambda
	r_n( B,f)\right]\mbox{d}\pi( B,f)} .$$
The notation ${\mbox{d}\hat{\rho}_{\lambda}}/{\mbox{d}\pi}$ means the density of $\hat{\rho}_{\lambda}$ with respect to $\pi$.

\subsection{Theoretical results}
As $\mathbb E[Y|X ]=f^{\star}( \left\langle X,B^{ \star } \right\rangle)$ almost surely, it is noted that for all $ ( B,f) \in \mathcal S_{1,+}^p \times \mathcal F_n(C+1) $,
\begin{align*}
R( B,f)- R ( B^\star,f^\star)
=
\mathbb E \left [Y-f( \left\langle X,B  \right\rangle)\right]^2-\mathbb E \left [Y-f^{\star}( \left\langle X,B^{\star } \right\rangle)\right]^2
= 
\mathbb E \left [ f( \left\langle X,B  \right\rangle)-f^{\star}( \left\langle X,B^{\star } \right\rangle)\right]^2
\end{align*}
(Pythagora's theorem).

For any positive integer $M\leq n$, we set
$$ 
\left(B^{\star}_{M},f^{\star}_{M}\right) 
\in
\arg\min_{( B,f)\in\mathcal{S}^{p}_{1,+} \times\mathcal{F}_{M}(C)}
R(B,f). 
$$
At this stage, it is very important to note that, for each $M$, the infimum $f^{\star}_{ M}$ is defined on $\mathcal F_{M}(C)$, whereas the prior charges a slightly bigger set, namely $\mathcal F_{M}(C+1)$.

The main result of the paper is the following theorem whose proof is given in Section \ref{sectionproofs}.. Here and everywhere, the
wording ``with probability $1-\delta$'' means the probability evaluated
with respect to the distribution $\mathbf P^{\otimes n}$ of the data {\it and} the conditional probability measure $\hat \rho_{\lambda}$.

\begin{theo}
	\label{thm2}
	Assume that Assumption $\mathbf N$ and Assumption $\mathbf B$ hold.  Set 
	$w = 64(C+1)\max[L,C+1], C_1 := 8 [(C+1)^2+ \sigma^{2}] $ and take 
	\begin{equation}
	\label{equationlambda}
	\lambda = \frac{n}{w + 2C_1}.
	\end{equation}
	Then, for all $\delta\in (0,1) $, with probability at least $1-\delta$	we have
	\begin{align*}
	R(\hat{B}_{\lambda},\hat{f}_{\lambda})- R(B^{\star},f^{\star}) 
	\leq
	\Xi\inf_{
		\tiny{\begin{array}{c} 
		1\leq M \leq n \end{array}}} \Bigg\{
	R(B^{\star}_{ M},f^{\star}_{ M}) -
	R(B^{\star},f^{\star}) 
	+
	\frac{ M\log(Cn) + {\rm rank}(B)d\log(16n) 
	+
	d\log d \log(2n e)
		+ 
	\log\left(\frac{2}{\delta}\right)}{n}
\Bigg\},
	\end{align*}
	where $\Xi$ is a positive constant, depends on $L$, $C$, $\sigma$ and $C_\phi $ only.
\end{theo}

As shown in \cite{alquier2013sparse}, when $ f^\star $ belongs to a So\-bo\-lev space, we can derive a more specific non-parametric rate for the above theorem. For example assume that
$\{\varphi_j\}_{j=1}^{\infty}$ is the (non-normalized) trigonometric system and that the target $f^{\star}$ be\-longs to the So\-bo\-lev ellipsoid, defined by
$$\mathcal W\left(k,\frac{{6}C^2}{\pi^2}\right)= \left\{ f\in L_{2}([-1,1]):
f=\sum_{j=1}^{\infty}\beta_j \varphi_j \mbox{ and } 
\sum_{j=1}^{\infty} j ^{2k}
\beta_{j}^{2}  \leq \frac{{6}C^2}{\pi^2} \right\}
$$
for some unknown regularity parameter $k\geq 2$ (see, e.g., Tsybakov
\cite{Tsybakov}). In this context, the approximation sets
$\mathcal F_M(C+1)$ take the form
\begin{align*}
 \mathcal F_M(C+1)
 = 
 \left\{ f\in L_{2}([-1,1]): f=\sum_{j=1}^{M}\beta_j
\varphi_j, \sum_{j=1}^{M} j|\beta_{j}|  \leq C+1 \mbox{ and } \beta_M\neq 0\right\}.
\end{align*}
As the regularity parameter $k$ is not given, the following results is in the so-called adaptive setting. However, we need to put the following additional assumption to obtain the result. The proof for Corollary \ref{thm3} is similar to Corollary 4 in \cite{alquier2013sparse} and thus omitted.

\noindent\textbf{Assumption $\mathbf D$.} The random variable $\left\langle X,B^\star \right\rangle $ has a probability density on
$[-1,1]$, bounded from above by a positive constant $A $.

\begin{cor}
	\label{thm3}
	Assume that Theorem \ref{thm2} and additional Assumption $\mathbf D$ hold. Suppose that $f^{\star}$ be\-lon\-gs to the So\-bo\-lev	ellipsoid $ \mathcal W(k,{6}C^2/\pi^2)$, where the unknown real number $k\geq 2$ is a
	regularity parameter. Take 
	$\lambda$ as in (\ref{equationlambda}).
	Then, for all $\delta\in (0,1) $, with probability at least $1-\delta$	we have
	\begin{align}
 R(\hat{B}_{\lambda},\hat{f}_{\lambda}) -
	R(B^{\star},f^{\star})\nonumber
\leq 
\Xi' \left\{
	\left(\frac{\log(Cn)}{n}\right)^{\frac{2k}{2k+1}} +
	\frac{{\rm rank}(B^\star)d\log(16n) 
		+
		d\log d \log(2n e)}{n}
+
\frac{	\log\left(\frac{2}{\delta}\right)}{n}
\right\},\label{kx}
	\end{align}
	where $\Xi'$ is a positive constant, function of $L$, $C$, $\sigma$,
	$C_\phi$ and $A $ only.
\end{cor}

From Theorem \ref{thm2}, it is actually possible to derive that the Gibbs posterior $ \hat{\rho}_\lambda $ contracts around $ (B^\star,f^\star) $ at the optimal rate.

\begin{cor}
	\label{thrm_contraction} 
	Under the same assumptions for Theorem~\ref{thm2}, and the same definition for $\tau$ and $\lambda^* $, let $\varepsilon_n $ be any sequence in $(0,1)$ such that $\varepsilon_n\rightarrow 0$ when $n\rightarrow\infty$. Define
	\begin{multline*}
	\mathcal{E}_n 
	= 
	\Biggl\{ ( B,f) \in \mathcal S_{1,+}^p \times \mathcal F_n(C+1) : 
	R( B,f) - R(B^{\star},f^{\star}) 
\leq
\Xi\inf_{
	\tiny{\begin{array}{c} 
		1\leq M \leq n \end{array}}} \Bigg\{
R(B^{\star}_{ M},f^{\star}_{ M}) -
R(B^{\star},f^{\star}) +
\\
\frac{ M\log(Cn) + {\rm rank}(B)d\log(16n) 
+
d\log d \log(2n e)
+ 
\log\left(\frac{2}{\varepsilon_n}\right)}{n}
\Bigg\},
	\end{multline*}
	Then
	$
	\mathbb{E} \Bigl[ \mathbb{P}_{ ( B,f) \sim \hat{\rho}_{\lambda}} ( ( B,f) \in\mathcal{E}_n) \Bigr] 
	\geq 
	1-\varepsilon_n \xrightarrow[n\rightarrow \infty]{} 1.  
	$ 
\end{cor}

\section{Proofs}
\label{sectionproofs}

For sake of simplicity in the proofs, we put
\begin{align*}
R^\star : = R ( B^\star,f^\star),
\quad
r_n^\star : =  r_n( B^\star,f^\star).
\end{align*}

\noindent Here we let $\pi$ be the prior probability measure on
$\mathbb R^p \times \mathcal F_n(C+1)$
equipped with its canonical Borel $\sigma$-field. We have that, for
each $f = \sum_{j=1}^M \beta_j \varphi_j \in \mathcal F_M(C+1)$,
\begin{align*}
\|f\|_{\infty} \leq \sum_{j=1}^M |\beta_j| \leq C+1.
\end{align*}

\noindent We start with some technical lemmas. Lemma \ref{lemmemassart} is a
version of Bernstein's inequality, whose proof can be found in 
\cite[inequality (2.21)]{Massart}. For a random variable $Z$, the
notation $(Z)_+$ means the positive part of $Z$.
\begin{lem}
	\label{lemmemassart}
	Let $T_{1}, \hdots, T_{n}$ be independent real-valued random variables. Assume
	that there exist two positive constants $v$ and $w$ such that, for all integers
	$k\geq 2$,
	$ \sum_{i=1}^{n} \mathbb{E}\left[(T_{i})_{+}^{k}\right] \leq
	\frac{k!}{2}vw^{k-2}.$
	Then, for any $\zeta\in (0,1/w) $,
	$$ \mathbb{E}
	\exp\left[\zeta\sum_{i=1}^{n}(T_{i}-\mathbb{E} T_{i}) \right]
	\leq \exp\left(\frac{v\zeta^{2}}{2(1-w\zeta)} \right) .$$
\end{lem}
Given a measurable space $(E,\mathcal E)$ and two probability measures $\mu_1$
and
$\mu_2$ on $(E,\mathcal E)$, we denote by $\mathcal K(\mu_1,\mu_2)$ the
Kullback-Leibler divergence of $\mu_1$ with respect to $\mu_2$, defined by
$
\mathcal K(\mu_1,\mu_2)
=
\int  \log \left (
\frac{\mbox{d}\mu_1}{\mbox{d}\mu_2}\right)\mbox{d}
\mu_1\, \mbox{if } \mu_1 \ll
\mu_2,$ 
and $
\infty $ otherwise. (Notation $\mu_1 \ll \mu_2$ means ``$\mu_1$ is absolutely continuous with respect to $\mu_2$''). Lemma \ref{lemmecatoni} is a classical result and its proof can be found for example in \cite[page 4]{catoni2007pac}.

\begin{lem}
	\label{lemmecatoni}
	Let $(E,\mathcal E)$ be a measurable space. For any probability measure $\mu$ on
	$(E,\mathcal E)$ and any measurable function
	$h:E\rightarrow\mathbb{R}$ such that $ \int (\exp \circ h)\emph{d}\mu 
	<\infty $,
	we have
	\begin{equation} 
	\log \int (\exp\circ h)\emph{d}\mu 
	=
	\sup_{m}\left(\int h\emph{d}m
	-
	\mathcal{K}(m,\mu)\right),
	\end{equation}
	where the supremum is taken over all probability measures on $(E,\mathcal E)$
	and, by convention,
	$\infty-\infty=-\infty$. Moreover, as soon as $h$ is bounded from above on the support of $\mu$, the supremum with respect to
	$m$ on the right-hand side of (\ref{lemmacatoni}) is reached for the Gibbs
	distribution $g$ given by
	$$\frac{\emph{d} g}{\emph{d} \mu}(e)=\frac{\exp\left[h(e)\right]}{\displaystyle
		\int(\exp \circ h)\emph{d}\mu},  \quad e \in E.$$
\end{lem}

\begin{lem}
	\label{thm1}
	Assume that Assumption {\bf N} holds. Set  $w= 16(C+1)\max[L,2(C+1)], C_1 := 8 [(C+1)^2+ \sigma^{2}] $ and take
	$
	\lambda \in\,\left(0,
	\frac{n}{w+ C_1 }\right)
	$ and put
	\begin{equation}
	\label{defalpha}
	\alpha = \left(\lambda
	-\frac{\lambda^{2} C_1 }{2n(1-\frac{ C_2 \lambda}{n})}\right)
	\quad \text{and} \quad
	\beta = \left(\lambda
	+\frac{\lambda^{2} C_1}{2n(1-\frac{ C_2 \lambda}{n})}\right) .
	\end{equation}
	Then, for any $\delta\in\,(0,1)$ and any data-dependent probability measure
	$\hat{\rho}_\lambda \ll\pi$ we have,
	\begin{align}
	\label{lemma1}
\mathbb{E}
 \int \exp\Biggl [
\alpha
\left(R( B,f)-R^\star \right)
+
\lambda\left(-r_n( B, f)+ r_n^\star \right)
- \log\left(\frac{\mbox{d}\hat{\rho}_\lambda}{\mbox{d}\pi}( B, f)\right)
- 
\log\frac{2}{\delta}
\Biggr
]\mbox{d}\hat{\rho}_\lambda( B,f)
\leq
\delta/2,
	\end{align}
	\begin{align}
	\label{lemma2}
\mathbb{E}
\sup_{\rho} \exp\Biggl [
\beta   \left( -\int R( B,f)\mbox{d}\rho -R^\star \right)
+
\lambda\left( \int r_n( B, f)\mbox{d}\rho - r_n^\star \right)
- 
\mathcal{K}(\rho,\pi)
- 
\log\frac{2}{\delta}
\Biggr]
\leq
\delta/2,
	\end{align}
	
\end{lem}

\noindent{\bf Proof of Lemma \ref{thm1}}. 
\\ 
Fix
$ B\in\mathcal{S}^{p}_{1,+}$ and $f \in \mathcal F_n(C+1)$. We start with an application of Lemma
\ref{lemmemassart} to the random variables
$$ T_{i} 
=  
- \left(Y_{i}-f(  \left\langle X, B \right\rangle )\right)^{2}
+ 
\left(Y_{i}-f^{\star}( \left\langle X, B^{\star} \right\rangle)\right)^{2}
, \quad i=1, \hdots, n.
$$
Note that $ T_{i},i=1, \hdots, n $ are independent and we have that
\begin{align*}
\sum_{i=1}^{n} \mathbb{E} T_{i}^{2}
&  = 
\sum_{i=1}^{n} \mathbb{E} \left \{
\left[2Y_{i} - f( \left\langle X, B \right\rangle)-f^{\star}( \left\langle X,B^\star \right\rangle)\right]^{2}
\left[f( \left\langle X,B \right\rangle)-f^{\star}( \left\langle X,B^{\star } \right\rangle)\right]^2
\right\}
\\
& =
\sum_{i=1}^{n} \mathbb{E} \left\{
\left[2\epsilon_{i} +f^{\star }( \left\langle X,B^{\star } \right\rangle) -
f( \left\langle X,B \right\rangle)\right]^{2}
\left[f( \left\langle X,B \right\rangle)-f^{\star}( \left\langle X,B^{\star } \right\rangle)\right]^2
\right\}
\\
& \leq
\sum_{i=1}^{n} \mathbb{E} \left\{
\left [ 8\epsilon_{i}^{2} + 8(C+1)^{2}\right]
\left[f( \left\langle X,B \right\rangle)-f^{\star}( \left\langle X,B^{\star } \right\rangle)\right]^2
\right\}.
\\
&\leq
8\left[(C+1)^2+ \sigma^2\right] \sum_{i=1}^{n}
\mathbb{E}\left[f( \left\langle X,B \right\rangle)-f^{\star}( \left\langle X,B^\star \right\rangle)\right]^2
:= v, 
\end{align*}
where we set
$
C_1 := 8 [(C+1)^2+ \sigma^{2}]; \,
\text{and} \,
v=  
nC_1
\left[R( B,f)-R^{\star}\right].
$
More generally, for all integers $k \geq 3$,
\begin{align*}
\sum_{i=1}^{n} \mathbb{E}\left[(T_{i})_{+}^{k}\right]
&  \leq 
\sum_{i=1}^{n} \mathbb{E} \left\{
\left|2Y_{i} - f( \left\langle X,B \right\rangle)-f^{\star}( \left\langle X,B^\star \right\rangle)\right|^{k}
\left|f( \left\langle X,B \right\rangle)-f^{\star}( \left\langle X,B^{\star } \right\rangle)\right|^k
\right\}
\\
& = \sum_{i=1}^{n} \mathbb{E} \left\{
\left |2\epsilon_{i} +f^{\star }( \left\langle X,B^{\star } \right\rangle) -
f( \left\langle X,B \right\rangle)\right|^{k}
\left|f( \left\langle X,B \right\rangle)-f^{\star}( \left\langle X,B^{\star } \right\rangle)\right|^k
\right\}\\
& \leq 2^{k-1} \sum_{i=1}^{n} \mathbb{E} \left\{
\left[2^k|\epsilon_{i}|^{k} + 2^k(C+1)^{k}\right]
2^{k-2}(C+1)^{k-2}
\left|f( \left\langle X,B \right\rangle)
-
f^{\star}( \left\langle X,B^* \right\rangle)\right|^{2}
\right\}.
\end{align*}
In the last inequality, we used the fact that $|a+b|^k\leq 2^{k-1}(|a|^k+|b|^k)$.
Therefore, by Assumption {\bf N},
\begin{align*}
\sum_{i=1}^{n} \mathbb{E}\left[(T_{i})_{+}^{k}\right]
& \leq 
\sum_{i=1}^{n} \left[2^{2k-2}{k!}\sigma^{2}L^{k-2}
+
2^{2k-1}(C+1)^{k}\right] 
2^{k-2}(C+1)^{k-2}
\left[R( B,f) - R^{\star}\right]
\\
& = 
v \times \frac{\left[2^{2k-2}{k!}\sigma^{2}L^{k-2}
	+
	2^{2k-1}(C+1)^{k}\right] 
	2^{k-2}(C+1)^{k-2}}
{[2(C+1)^2+4\sigma^{2}]}
\\
& \leq v \times
\frac{ k! 8^{k-2} \max\left [L^{k-2},
	2^{k-2}(C+1)^{k-2}\right]
	2^{k-2}(C+1)^{k-2}}{2}
:= \frac{k!}{2}v w^{k-2},
\end{align*}
with $w= 64(C+1)\max[L,C+1] $. Thus, for any $\lambda \in (0,n/w ) $, taking
$\zeta=\lambda/n$, we apply Lemma \ref{lemmemassart} to get
\begin{align*}
\mathbb{E} \exp\left[ \lambda
\left(R( B,f)- R^{\star}
-
r_n( B,f) + r_n^{\star}\right)\right]
\leq
\exp\left(\frac{v\lambda^{2}}{2n^{2}(1-\frac{w\lambda}{n})}\right)
=
\exp\left(\frac{C_1
\left[R( B,f)-R^{\star}\right]\lambda^{2}  }{ 2n (1-\frac{w\lambda}{n})}\right)
.
\end{align*}
Therefore, we obtain, with $ \alpha $ given in \eqref{defalpha},
\begin{align*}
\mathbb{E} \exp \Bigg [
\alpha
\left(R( B,f)- R^\star \right)
+\lambda\left(-r_n( B,f)+ r_n^\star \right)
- \log\left(\frac{2}{\delta}\right)
\Bigg]
\leq \delta/2.
\end{align*}
Next, integrating with respect to $ \pi $
and consequently using Fubini's theorem, we obtain
\begin{align*}
\mathbb{E} \int \exp\Biggl[
\alpha
\left(R( B,f)-R^\star \right)   
+
\lambda\left(-r_n( B,f)+ r_n^\star \right)-
\log\left(\frac{2}{\delta}\right)\Biggr] \mbox{d}\pi( B,f)
\leq \delta/2.
\end{align*}
To obtain \eqref{lemma1}, it is noted that for any measurable function $ h $,
\begin{align*}
\int \exp [ h(B,f) ] d \pi
 = 
 \int \exp \left[ h(B,f) - \log \frac{d \hat{\rho}_\lambda }{d \pi} (B,f) \right] d \hat{\rho}_\lambda .
\end{align*}

The proof for \eqref{lemma2} is similar. More precisely, we apply Lemma \ref{lemmemassart} with $T_i=  (Y_{i}-f( \left\langle X,B \right\rangle))^{2}
- (Y_{i}-f^{\star}( \left\langle X,B^{\star } \right\rangle))^{2} $. We obtain, for any $\lambda \in (0,n/w ) $,
\begin{align*}
\mathbb{E} \exp\left[ \lambda
\left(r_n( B,f) + r_n^{\star} 
-
R( B,f) + R^{\star}
\right)\right]
\leq
\exp\left(\frac{v\lambda^{2}}{2n^{2}(1-\frac{w\lambda}{n})}\right).
\end{align*}
By rearranging terms, using definition of $ \beta $ in \eqref{defalpha}, and multiplying both side by $ \delta/2 $, we obtain
\begin{align*}
\mathbb{E} \exp \Bigg [
\beta
\left( - R( B,f) + R^\star \right)
+
\lambda\left(r_n( B,f) - r_n^\star \right)
- \log \frac{2}{\delta} 
\Bigg]
\leq \delta/2.
\end{align*}
Integrating with respect to $ \pi $ and using Fubini's theorem to get,
\begin{align*}
\mathbb{E} \int \exp \Bigg [
\beta
\left( - R( B,f) + R^\star \right)
+
\lambda\left(r_n( B,f) - r_n^\star \right)
- \log \frac{2}{\delta} 
\Bigg]  \mbox{d}\pi
\leq 
\delta/2.
\end{align*}
Now, Lemma \ref{lemmecatoni} is applied on the integral and this directly yields \eqref{lemma2}.
\hfill $\blacksquare$

\begin{proof}[\textbf{Proof of Theorem \ref{thm2}}]
Remind that $\mathbf P^{\otimes n} $ stands for the distribution of the sample $ \mathcal{D}_n $, the equation \eqref{lemma1} can be written conveniently as
\begin{align*}
\mathbb{E}_{\mathcal{D}_n \sim \mathbf P^{\otimes n} }
\mathbb{E}_{ (\hat{ B},\hat{f}) \sim \hat{\rho}_{\lambda}} \exp\Biggl [
\alpha
\left(R( \hat{ B},\hat{f})-R^\star \right)
+
\lambda\left(-r_n( \hat{ B},\hat{f})
+ r_n^\star \right)
- \log\left(\frac{\mbox{d}\hat{\rho}_\lambda}{\mbox{d}\pi}( \hat{ B},\hat{f})\right)
- 
\log\frac{2}{\delta}
\Biggr] 
\leq
\delta/2,
\end{align*}

Now, we use the standard Chernoff's trick to transform an exponential moment inequality into a deviation inequality, that is using the inequality $\exp(\lambda x) \geq
\mathbf{1}_{\mathbb{R}_{+}}(x)$. We obtain, with probability at least $1-\delta/2 $, for any $\delta\in\,(0,1)$ and any data-dependent probability measure
$\hat{\rho}_\lambda \ll\pi$,
\begin{align*}
 R( \hat{ B},\hat{f})  -R^\star  
\leq 
\frac{\lambda}{\alpha} 
\Biggr(
 r_n( \hat{ B},\hat{f}) 
 - r_n^\star 
+\frac{
\log\left(\frac{\mbox{d}\hat{\rho}_\lambda}{\mbox{d}\pi}( \hat{ B},\hat{f})\right)
	+
	\log\left(\frac{2}{\delta}\right)}{\lambda}\Biggr).
\end{align*}
It is noted that we have
\begin{align*}
\log\left(\frac{\mbox{d}\hat{\rho}_\lambda}{\mbox{d}\pi}( \hat{ B},\hat{f})\right)
=
\log\left(\frac{\exp (-\lambda r_n(\hat{ B},\hat{f}) ) }{ \int \exp (-\lambda r_n(B,f) ) \mbox{d}\pi } \right)
= 
-\lambda r_n(\hat{ B},\hat{f})
-
\log \int \exp (-\lambda r_n(B,f) ) \mbox{d}\pi 
,
\end{align*}
thus we obtain, with probability at least $1-\delta/2 $,
\begin{align*}
R( \hat{ B},\hat{f}) 
-
R^\star  
\leq 
\frac{1}{\alpha} 
\Biggr(
\log \int \exp (-\lambda r_n(B,f) ) \mbox{d}\pi 
- \lambda r_n^\star 
+
\log\left(\frac{2}{\delta}\right)
\Biggr).
\end{align*}
Next, using Lemma \ref{lemmecatoni} we deduce that, with probability at least $1-\delta/2 $ 
\begin{align}
\label{unionbound1}
R( \hat{ B},\hat{f}) 
-
R^\star  
\leq 
\frac{\lambda}{\alpha} 
\Biggr( \int r_n(B,f)  \mbox{d}\hat{\rho}_\lambda
- r_n^\star 
+\frac{
\mathcal{K}(\hat{\rho}_\lambda,\pi)
	+
	\log\left(\frac{2}{\delta}\right)}{\lambda}\Biggr).
\end{align}

Now, from \eqref{lemma2}, an application of  the standard Chernoff's trick, that is using the inequality $\exp(\lambda x) \geq
\mathbf{1}_{\mathbb{R}_{+}}(x)$. We obtain, with probability at least $1-\delta/2 $, for any $\delta\in\,(0,1)$ and any data-dependent probability measure
$\hat{\rho}_\lambda \ll\pi$,
\begin{align}
\label{unionbound2}
\int r_n( B, f) \mbox{d}\hat{\rho}_\lambda
- r_n^\star 
\leq 
\frac{\beta}{\lambda} 
\Biggr(
\int R( B,f) \mbox{d}\hat{\rho}_\lambda 
-R^\star   \Biggr)
+\frac{
\mathcal{K}(\hat{\rho}_\lambda,\pi)
	+
\log\left(\frac{2}{\delta}\right)}{\lambda}
.
\end{align}
Combining \eqref{unionbound1} and \eqref{unionbound2} with a union bound argument gives the bound, with
probability at least $1-\delta $,
\begin{align}
  R( \hat{ B},\hat{f}) - R^\star  
\leq 
\inf_{ \rho }\Bigg \{ 
\frac{\beta}{\alpha} 
\bigg(\int R( B,f)
\mbox{d} \rho
-
R^{\star} \bigg)
+ 
2 \,\frac{
\mathcal{K}( \rho ,\pi)
	+ \log\left(\frac{2}{\delta}\right)
}{\alpha}
\Bigg\}.
\label{inegalitedegueu}
\end{align}
The last steps of the proof consist in deriving the right-hand side in the above inequality more explicit. In order to do that, we restrict the infimum bound to the following specific distribution.

Put $ B^{\star}_{ M} = U\Lambda U^\top $ and let $r = \#\{ i : \Lambda_i > \varepsilon \},$
with small $\varepsilon \in (0,1)$. Take
$$
\mbox{d}\rho^{1}_{\eta}
 \propto 
 \mathbf{1} (
\forall i: |v_i - \Lambda_i| \leq \varepsilon; 
\forall i = 1,\dots,r: \| u_i -
U_i \|_F \leq \eta) \pi (du,dv)
$$
For any positive integer $M\leq n$ and any
$\eta,\gamma\in\,]0,1/{n}]$, let the probability measure
$\rho_{M,\eta,\gamma}$ 
be defined by
$$ 
\mbox{d}\rho_{ M,\eta,\gamma}( B,f) 
=
\mbox{d}\rho^{1}_{\eta}( B) \mbox{d}\rho^{2}_{ M,\gamma}(f),
$$
with
$$
\frac{\mbox{d}\rho^{2}_{ M,\gamma}}{ \mbox{d}\nu_{M}}(f)  \propto
\mathbf{1}_{[\|f-f^{\star}_{ M}\|_M \leq \gamma]} $$
where, for $f=\sum_{j=1}^{M}\beta_{j}\varphi_{j}\in\mathcal{F}_M(C+1)$, we put
$$ \|f\|_M = \sum_{j=1}^{M}j|\beta_{j}| .$$
With this notation, inequality \eqref{inegalitedegueu} leads to
\begin{align}
  R( \hat{ B},\hat{f}) - R^\star  
\leq 
\inf_{1\leq M \leq n }
\inf_{\eta,\gamma>0} 
\Bigg\{        
\frac{\beta}{\alpha} 
\Bigg(\int
R( B,f)\mbox{d}\rho_{ M,\eta,\gamma}( B,f) 
-
R^{\star} \Bigg)
 +
2\,\frac{\mathcal{K}(\rho_{ M,\eta,\gamma}, \pi)
	+ 
	\log\left(\frac{2}{\delta}\right)}{\alpha}
\Biggr\}.
\label{inegalitedegueu2}
\end{align}
To finish the proof, we have to control the different terms in
\eqref{inegalitedegueu2}. Note first that
\begin{align*}
\mathcal{K}(\rho_{ M,\eta,\gamma}, \pi)
=
\mathcal{K}(\rho^{1}_{ \eta}\otimes\rho^{2}_{ M,\gamma},\mu \otimes\nu_{M}) 
= 
\mathcal{K}(\rho^{1}_{\eta},\mu) 
+
\mathcal{K}(\rho^{2}_{ M,\gamma},\nu_{M})
+
 \log \frac{1 - (1/10)^n}{10^{-M}}
\nonumber .
 \end{align*}
By technical Lemma \ref{technique1}, we know that 
$$
\mathcal{K}(\rho^{1}_{\eta},\mu)
\leq
	 rd\log(16/\eta) +
C_{D_1}d\log d(1 + \log(2/\varepsilon))
$$
By technical Lemma 10 in \cite{alquier2013sparse}, we have
$$ 
\mathcal{K}(\rho^{2}_{ M,\gamma},\nu_{M})  
= 
M \log\left(\frac{C+1}{\gamma}\right).
$$
Putting all the pieces together, we are led to
\begin{equation}
\mathcal{K}(\rho_{ M,\eta,\gamma},\pi)
\leq 
 rd\log(1/c) 
 +
C_{D_1}d\log d(1 + \log(2/\delta))
+ 
M \log \left(\frac{C+1}{\gamma}\right)
+
\log \frac{1}{10^{-M}}
. 
\label{L2}
\end{equation}
Finally, it remains to control the term
$\int R( B,f)\mbox{d}\rho_{ M,\eta,\gamma}( B,f).$
To this aim, we write
\begin{align*}
\int R( B,f)\mbox{d}\rho_{ M,\eta,\gamma}( B,f)
& =
\int\mathbb{E} \left[\left(Y-f( \left\langle X,B  \right\rangle)\right)^{2}\right]
\mbox{d}\rho_{ M,\eta,\gamma}( B,f)
\\
& =
\int\mathbb{E} \big[\big(Y-f^{\star}_{ M}( \left\langle X,B^{ \star }_{ M}  \right\rangle )+
f^{\star}_{ M}( \left\langle X,B^{ \star }_{ M}  \right\rangle )-f( \left\langle X,B^{ \star }_{ M}  \right\rangle )
+
f( \left\langle X,B^{ \star }_{ M}  \right\rangle )-f( \left\langle X,B  \right\rangle)\big)^{2}\big]
\mbox{d}\rho_{ M,\eta,\gamma}( B,f)\\
& =
R( B^{\star}_{ M},f^{\star}_{ M}) + \int\mathbb{E} \Bigl[
\left(f^{\star}_{ M}( \left\langle X,B^{ \star }_{ M}  \right\rangle )- 
f( \left\langle X,B^{ \star }_{ M}  \right\rangle )\right)^{2}
+
 \left(f( \left\langle X,B^{ \star }_{ M}  \right\rangle )-f( \left\langle X,B  \right\rangle)\right)^{2}
\\
&  +
 2 \left(Y-f^{\star}_{ M}( \left\langle X,B^{ \star }_{ M}  \right\rangle )\right)
\left(f^{\star}_{ M}(
 \left\langle X,B^{ \star }_{ M}  \right\rangle )-f(  \left\langle X,B^{ \star }_{ M}  \right\rangle )\right)
 + 
  2\left(Y-f^{\star}_{ M}(  \left\langle X,B^{ \star }_{ M}  \right\rangle )\right)\left(f( \left\langle X,B^{ \star }_{ M}  \right\rangle )-f( \left\langle X,B  \right\rangle)\right)
\\
& + 2 \left(f^{\star}_{ M}( \left\langle X,B^{ \star }_{ M}  \right\rangle )
-
f(  \left\langle X,B^{ \star }_{ M}  \right\rangle )\right)\left(f( \left\langle X,B^{ \star }_{ M}  \right\rangle )
-
f( \left\langle X,B  \right\rangle)\right)
\Bigr] \mbox{d}\rho_{ M,\eta,\gamma}( B,f)
\\
& := R( B^{\star}_{ M},f^{\star}_{ M}) + \mathbf A+\mathbf B+\mathbf
C+\mathbf D+\mathbf E.
\end{align*}

\noindent Computation of C:
By Fubini's theorem,
\begin{align*}
{\bf C}  &= \mathbb{E}\left[\int 2 \left(Y-f^{\star}_{ M}(  \left\langle X,B^{ \star }_{ M}  \right\rangle )\right)
\left(f^{\star}_{ M}( \left\langle X,B^{ \star }_{ M}  \right\rangle )
-
f(  \left\langle X,B^{ \star }_{ M}  \right\rangle )\right) \mbox{d}\rho_{ M,\eta,\gamma}( B,f) \right]
\\
&= 
\mathbb{E}\Biggl\{\int \Biggl[2 \left(Y-f^{\star}_{ M}(  \left\langle X,B^{ \star }_{ M}  \right\rangle )\right)
\int
\left(f^{\star}_{ M}( \left\langle X,B^{ \star }_{ M}  \right\rangle )
-
f(  \left\langle X,B^{ \star }_{ M}  \right\rangle )\right)
\mbox{d}\rho^{2}_{ M,\gamma}(f)
\Biggr]
\mbox{d}\rho^{1}_{ \eta}( B)
\Biggr\}.
\end{align*}
By the triangle inequality, for 
$
f=\sum_{j=1}^M \beta_j \varphi_j
$ 
and 
$
f^{\star}_{ M}=\sum_{j=1}^M (\beta^{\star}_{ M})_j \varphi_j
$, it holds
$
\sum_{j=1}^M j |\beta_j|\leq \sum_{j=1}^M j \left | \beta_j-(\beta^{\star}_{ M})_j \right|+ \sum_{j=1}^M j \left |(\beta^{\star}_{ M})_j\right|.
$
Since $f^{\star}_{ M} \in \mathcal F_M(C)$, we have $\sum_{j=1}^M j | (\beta^{\star}_{ M})_j| \leq C$,
so that $\sum_{j=1}^M j |\beta_j|\leq C+1$ as soon as $\| f-f^{\star}_{ M}\|_M \leq 1$. This shows that the set
$$\left\{f = \sum_{j=1}^M \beta_j \varphi_j : \| f-f^{\star}_{ M}\|_M \leq \gamma\right\}$$
is contained in the support of $\nu_M$. In particular, this implies that $\rho^{2}_{ M,\gamma}$ is centered at $f^{\star}_{ M}$ and, consequently, 
$$
\int 
\left(f^{\star}_{ M}( \left\langle X,B^{ \star }_{ M}  \right\rangle ) 
-
f( \left\langle X,B^{ \star }_{ M}  \right\rangle )\right)
\mbox{d}\rho^{2}_{ M,\gamma}(f) = 0.
$$
This proves that $\mathbf C=0 $.

\noindent Control of A:
Clearly,
$
{\bf A} \leq  \int \sup_{y\in\mathbb{R}} \left((f^{\star}_{ M}(y)-f(y)\right)^{2}
\mbox{d}\rho^{2}_{ M,\gamma}(f)
\leq \gamma^{2}.
$

\noindent Control of B:
We have
\begin{align*}
{\bf B} 
= \int\mathbb{E} &\left[\left(f( \left\langle X,B^{ \star }_{ M}  \right\rangle )-f( \left\langle X,B   \right\rangle)\right)^{2} \right]
\mbox{d}\rho_{ M,\eta,\gamma}( B,f)
\leq 
\int \mathbb{E} \left[\left(C_\phi (C+1) 
( B^{\star }_{ M}- B )X \right)^2\right] \mbox{d}\rho^{1}_{ \eta}( B)
\quad (\mbox{by the mean value theorem})
\\
&\leq 
C_\phi^2(C+1)^{2}\mathbb{E}\left [\|X \|^2_{\infty}\right] 
\int
\| B^{\star}_{ M}- B\|_{F}^{2}
 \mbox{d}\rho^{1}_{\eta}( B)
(\mbox{by Assumption } \mathbf D).
\end{align*}
Using Lemma 6 from \cite{mai2017pseudo}, we have that $  \int
\| B^{\star}_{ M}- B\|_{F}^{2}
\mbox{d}\rho^{1}_{\eta}( B) 
\leq (3dc + 2r\eta )^2  $. Thus $ {\bf B} \leq C_\phi^2(C+1)^{2}(3dc + 2r\eta )^2 $.

\noindent Control of E:
Write
\begin{align*}
|{\bf E}| 
&\leq 
2 \int\mathbb{E}
\Bigl[
\left|f^{\star}_{ M}( \left\langle X,B^{ \star }_{ M}  \right\rangle )
-
f( \left\langle X,B^{ \star }_{ M}  \right\rangle )\right|
\left|f( \left\langle X,B^{ \star }_{ M}  \right\rangle )
-
f( \left\langle X,B  \right\rangle)\right|
\Bigr]
\mbox{d}\rho_{ M,\eta,\gamma}( B,f)
\\
&\leq
2 \int\mathbb{E}
\Bigl[
\left|f^{\star}_{ M}( \left\langle X,B^{ \star }_{ M}  \right\rangle )
-
f( \left\langle X,B^{ \star }_{ M}  \right\rangle )\right|
C_\phi(C+1)\left|( B^{\star}_{ M}- B )X \right|
\Bigr]
\mbox{d}\rho_{ M,\eta,\gamma}( B,f)
\\
&\leq 
2 \left(\int\mathbb{E} \left(f^{\star}_{ M}
(  \left\langle X,B^{ \star }_{ M}  \right\rangle )
-
f( \left\langle X,B^{ \star }_{ M}  \right\rangle )\right)^{2}
\mbox{d}\rho_{ M,\eta,\gamma}( B,f)\right)^{\frac{1}{2}}
\left(\int\mathbb{E} 
\left(C_\phi(C+1) ( B^{\star
	}_{ M}- B )X \right)^{2}
\mbox{d}\rho_{ M,\eta,\gamma}( B,f)\right)^{\frac{1}{2}}
\\
&\quad 
(\mbox{by the Cauchy-Schwarz inequality})\\
&
\leq 
2 \left(\gamma^{2}\right)^{\frac{1}{2}}
\left(C_\phi^2(C+1)^{2}(3dc + 2r\eta )^2\right)^{\frac{1}{2}}
= 
2 C_\phi(C+1) \gamma(3d \varepsilon + 2r\eta ).
\end{align*}

\noindent Control of D:
Finally,
\begin{align*}
{\bf D} 
&= 
2 \int\mathbb{E}
\left[
\left(Y-f^{\star}_{ M}( \left\langle X,B^{ \star }_{ M}  \right\rangle )\right)
\left(f( \left\langle X,B^{ \star }_{ M}  \right\rangle )-f( \left\langle X,B  \right\rangle)\right)
\right]
\mbox{d}\rho_{ M,\eta,\gamma}( B,f)
\\
&= 2 \int\mathbb{E}
\left[
\left(Y-f^{\star}_{ M}( \left\langle X,B^{ \star }_{ M}  \right\rangle )\right)\left(f^{\star}_{ M}( \left\langle X,B^{ \star }_{ M}  \right\rangle )
-f^{\star}_{ M}( \left\langle X,B  \right\rangle)\right)
\right]
\mbox{d}\rho^{1}_{ \eta}( B)
 \quad (\mbox{since }\int f\mbox{d}\rho^{2}_{ M,\gamma}(f)=f^{\star}_{ M})\\
&= 
2 \mathbb{E}
\left[
\left(Y-f^{\star}_{ M}( \left\langle X,B^{ \star }_{ M}  \right\rangle )\right)\int\left(f^{\star}_{ M}( \left\langle X,B^{ \star }_{ M}  \right\rangle )-f^{\star}_{ M}( \left\langle X,B  \right\rangle)\right)
\mbox{d}\rho^{1}_{ \eta}( B)
\right]
\\
&\leq
2 \sqrt{\mathbb{E}\left[
	\left(Y-f^{\star}_{ M}( \left\langle X,B^{ \star }_{ M}  \right\rangle )\right)^{2}\right]}
\sqrt{\mathbb{E}\left[
	\int\left(f^{\star}_{ M}( \left\langle X,B^{ \star }_{ M}  \right\rangle )-f^{\star}_{ M}( \left\langle X,B  \right\rangle)\right)
	\mbox{d}\rho^{1}_{ \eta}( B)
	\right]^{2}}
 \\
&= 2 \sqrt{R( B^{\star}_{ M},f^{\star}_{ M})} \sqrt{\mathbb{E}\left[
	\int\left(f^{\star}_{ M}( \left\langle X,B^{ \star }_{ M}  \right\rangle )-f^{\star}_{ M}( \left\langle X,B  \right\rangle)\right)
	\mbox{d}\rho^{1}_{\eta}( B)
	\right]^{2}}.
\end{align*}
As we have that
$
\left | f^{\star}_{ M}( \left\langle X,B^{ \star }_{ M}  \right\rangle )-f^{\star}_{ M}( \left\langle X,B  \right\rangle) \right|
\leq 
C_\phi (C+1) \left |\left\langle( B^{ \star }_{ M}- B )X\right\rangle \right|
\leq 
C_\phi (C+1) \| B^{\star}_{ M}- B\|_F,
$
it leads to
\begin{align*}
\left [\int \left(f^{\star}_{ M}(\left\langle X,B^{ \star }_{ M}  \right\rangle  )-f^{\star}_{ M}( \left\langle X,B  \right\rangle)\right)
\mbox{d}\rho^{1}_{ \eta}(B)\right]^2
\leq 
C_\phi ^2 (C+1)^2 \left [ \int  \| B^{\star}_{ M}- B\|_F \mbox{d}\rho^{1}_{\eta}( B)\right]^2
\leq 
C_\phi^2(C+1)^2 (3dc + 2r\eta )^2,
\end{align*}
and therefore
\begin{align*}
\mathbf D 
\leq 
2 C_\phi (C+1) (3dc + 2r\eta ) \sqrt{R(0,0)/2}
\leq 
\sqrt 2 C_\phi (C+1)  (3d \varepsilon + 2r\eta )\sqrt{C^{2}+\sigma^{2}}.
\end{align*}
Thus, taking $\eta= \gamma= \varepsilon= 1/n$ and putting all the pieces together, we obtain
$$ \mathbf A+\mathbf B+\mathbf C+\mathbf D+\mathbf E \leq \frac{\Xi_{1}}{n},$$
where
$ \Xi_{1}$ is a positive constant, function of $C$, $\sigma$ and $C_\phi$. Combining this inequality with \eqref{inegalitedegueu2}-\eqref{L2} yields,
with probability larger than $1-\delta $,
\begin{align*}
R(\hat B_{\lambda},\hat f_{\lambda})-R^\star 
   \leq 
\inf_{1\leq M \leq n }
\Biggl\{ 
\frac{\beta}{\alpha} 
\Bigg( R( B^{\star}_{ M},f^{\star}_{ M})
  -
  R^{\star} + \frac{
	\Xi_{1}}{n}\Bigg) 
+
2\,\frac{ M\log(10(C+1)n)	
+
rd\log(16n) +	C_{D_1}d\log d \log(2ne)
	+ \log\left(\frac{2}{\delta}\right)}{\lambda}
\Biggr\}.
\end{align*}
Choosing finally 
$$\lambda = \frac{n}{w+2	C_1 },$$
we obtain that there exists a positive constant $\Xi_2$, function of $L$,
$C$, $\sigma$ and $C_\phi$ such that,
with probability at least $1-\delta$,
\begin{align*}
R(\hat B_{\lambda},\hat f_{\lambda})-R^\star       
\leq 
\Xi_2
\inf_{  1\leq M \leq n }
\Biggl\{
R( B^{\star}_{ M},f^{\star}_{ M}) 
-
R^{\star}
  + 
  \frac{ M\log(10Cn)
	+  
rd\log(16n) +\Xi_3 d\log d \log(2ne)
	+ 
	\log\left(\frac{2}{\delta}\right)}{n}
\Biggr\}.
\end{align*}
This concludes the proof of Theorem \ref{thm2}.
\end{proof}

\begin{lem}
	\label{technique1}
Let $r = \#\{ i : \Lambda_i > \varepsilon \},$
with small $\varepsilon \in [0,1)$. Take
$$
\mbox{d}\rho^{1}_{\eta}
\propto 
\mathbf{1} (
\forall i: |v_i - \Lambda_i| \leq \varepsilon; 
\forall i = 1,\dots,r: \| u_i -
U_i \|_F \leq \eta) \mu (du,dv)
$$
	Then
	$$
	\mathcal{K}(\rho^{1}_{\eta},\mu)
	\leq
	 rd\log(16/\eta) +
	C_{D_1}d\log d \log(2e/\varepsilon)
	$$
 where $ C_{D_1} $ is a universal constant and depends only on $ D_1 $.
\end{lem}

\noindent{\bf Proof}.
We have that
\begin{align*}
\mathcal{K}(\rho^{1}_{ \eta},\mu)
 = &
  \log\frac{1}{\mu(\{ u,v: \forall i: |v_i - \Lambda_i| \leq \varepsilon;
  	 \forall i = 1,r: \| u_i - U_i \|_F \leq \eta \})}
\\
= & \log\frac{1}{ \mu \left( \left\lbrace
	\forall i = 1,r : \|u_{i.} - U_{i.}\|_{F}
	\leq \eta \right\rbrace \right)}
+
 \log\frac{1}{\mu(\{\forall i : |v_i - \Lambda_i |
	\leq \varepsilon\})}
 .
\end{align*}
The first log term
\begin{align*}
\pi \left( \left\lbrace \forall i = 1,r :
\|u_{i.} - U_{i.}\|_{F} \leq  \eta \right\rbrace \right)
\geq \prod_{i =1}^r  \Bigg[  \dfrac{\pi^{(d-1)/2}
	(\eta/2)^{d-1} }{ \Gamma(\frac{d-1}{2}+1)}     \Bigg/
\dfrac{2 \pi^{(d+1)/2}}{\Gamma(\frac{d+1}{2})}  \Bigg]  
 \geq  
\Bigg[ \dfrac{\eta^{d-1}}{2^d\pi}  \Bigg]^r
\geq   
\dfrac{\eta^{r(d-1)}}{2^{4rd}}   .
\end{align*}
Note for the above calculation: firstly the distribution of orthogonal vector is approximate by the uniform distribution on the sphere \cite{goldstein2017any} and secondly that the probability is greater or equal to the
volume of the (d-1)-"circle" with radius $ c/2 $
over the surface area of the $d$-``unit-sphere".

It is noted that if $ \gamma \sim Beta(a, b ) $ then $ \gamma^{1/2}$ has the pdf as 
$
 f(\gamma) 
 =
  2\frac{\gamma^{ 2 a -1} (1 - \gamma^2)^{b -1}}{Be( a, b )}, 0<\gamma<1 
 $ 
 where $ Be( a, b ) $ is the beta function. The second log term in the Kullback-Leibler term with $ a = \alpha_i, b = \sum_{i=1}^{d}\alpha_i - \alpha_i, \alpha_i = 1/d $
\begin{align*}
\pi(\{\forall i : |v_i - \Lambda_i | \leq \varepsilon \})
 =
\prod_{i=1}^d
\int_{\max(\Lambda_i - \varepsilon ,0)}^{\min(\Lambda_i + \varepsilon ,1)} 
 \frac{v_i^{ 2 a -1} (1 - v_i^2)^{b -1}}{ 2 Be( a, b )}  dv_i
\geq
\prod_{i=1}^d
\int_{0}^{\varepsilon} 
\frac{v_i^{ 2 a -1} (1 - v_i^2)^{b -1}}{ 2 Be( a, b )}  dv_i
\geq
\Xi_3 (\varepsilon /2d )^d \  e^{-d \log d}.
\end{align*}
 The interval of integration contains at least
an interval of length $ \varepsilon $.
Thus, we obtain
\begin{align*}
\mathcal{K}(\rho^{1}_{ \eta},\mu)
 \leq    
  \log \dfrac{2^{4rd}}{\eta^{r(d-1)}}
+ 
 \log\left( \frac{(2d)^d e^{d \log d }}{\Xi_3 \varepsilon^d}  \right)
 \leq  
  rd\log(\frac{16}{\eta}) +
\Xi_3d\log d \log(\frac{e2}{\varepsilon})
\end{align*}
for some absolute numerical constant $\Xi_3 $ that does not depend on $ r,n $ nor $ d $.
\hfill $\blacksquare$

\begin{proof}[\textbf{Proof of Corollary \ref{thrm_contraction}}]
We also start with an application of Lemma~\ref{thm1}, and focus on~\eqref{lemma1}, applied to $\delta :=\varepsilon_n $, that is:
	\begin{align*}
\mathbb{E}
\int \exp\Biggl [
\alpha
\left(R (B,f) - R^\star \right)
+
\lambda\left(-r_n( B, f)+ r_n^\star \right)
- \log\left(\frac{\mbox{d}\hat{\rho}_\lambda}{\mbox{d}\pi}( B, f)\right)
- 
\log\frac{2}{\varepsilon_n}
\Biggr
]\mbox{d}\hat{\rho}_\lambda( B,f)
\leq
\varepsilon_n/2
	\end{align*}
Using Chernoff's trick, this gives:
	$$
	\mathbb{E} \Bigl[ \mathbb{P}_{ (B,f) \sim \hat{\rho}_{\lambda}} ( (B,f) \in\mathcal{A}_n) \Bigr]
	\geq 1-\frac{\varepsilon_n}{2}
	$$
	where
	$$
	\mathcal{A}_n = \left\{ (B,f) : \alpha    \Bigl( R (B,f) - R^\star \Bigr)
	+\lambda\Bigl( -r_n(B,f) + r_n^\star \Bigr)      \leq      \log \left[\frac{d\hat{\rho}_{\lambda}}{d \pi} (B,f)  \right]
	+ \log\frac{2}{\varepsilon_n} \right\}.
	$$
Using the definition of $\hat{\rho}_\lambda $, for $ (B,f) \in \mathcal{A}_n $ we have
	\begin{align*}
	\alpha    \Bigl( R(B,f) - R^\star \Bigr)
	&    \leq  
	\lambda\Bigl( r_n(B,f) - r_n^\star \Bigr)  +       \log \left[\frac{d\hat{\rho}_{\lambda}}{d \pi} (B,f)  \right]
	+ \log\frac{2}{\varepsilon_n}
	\\
	& \leq -\log\int\exp\left[-\lambda r_n (B,f) \right]\pi({\rm d}(B,f)) - \lambda r_n^\star
	+ \log\frac{2}{\varepsilon_n}
	\\
	& = \lambda\Bigl( \int r_n(B,f) \hat{\rho}_{\lambda}({\rm d}(B,f)) - r_n^\star \Bigr)  +    \mathcal{K}(\hat{\rho}_\lambda,\pi)
	+ \log\frac{2}{\varepsilon_n}
	\\
& = 
	\inf_{\rho} \left\{ \lambda\Bigl( \int r_n(B,f) \rho({\rm d}(B,f)) - r_n^\star \Bigr)  +    \mathcal{K}(\rho,\pi)
	+ \log\frac{2}{\varepsilon_n} \right\}.
	\end{align*}
	
\noindent Now, let us define
	$
	\mathcal{B}_n 
	:= 
	\left\{\forall\rho\text{, }\beta \left(-\int R(B,f) d\rho + R^\star \right)
	+ \lambda \left( \int r_n d\rho - r_n^\star \right) \leq
	\mathcal{K}(\rho, \pi) + \log \frac{2}{\varepsilon_n}\right\}. $
	
	\noindent Using~\eqref{lemma2}, we have that
	$$
	\mathbb{E} \Bigl[\mathbf{1}_{\mathcal{B}_n} \Bigr]
	\geq 1-\frac{\varepsilon_n}{2}.
	$$
	We will now prove that, if $\lambda$ is such that $\alpha>0 $,
	$$
	\mathbb{E} \Bigl[ \mathbb{P}_{ (B,f) \sim \hat{\rho}_{\lambda}} ( (B,f) \in\mathcal{E}_n) \Bigr] 
	\geq
 \mathbb{E} \Bigl[ \mathbb{P}_{ (B,f) \sim \hat{\rho}_{\lambda}} ( (B,f) \in\mathcal{A}_n)\mathbf{1}_{\mathcal{B}_n} \Bigr]
	$$
	which, together with
	\begin{align*}
	\mathbb{E} \Bigl[ \mathbb{P}_{ (B,f)\sim \hat{\rho}_{\lambda}} ( (B,f) \in\mathcal{A}_n) \mathbf{1}_{\mathcal{B}_n} \Bigr]
	 = \mathbb{E} \Bigl[ (1-\mathbb{P}_{ (B,f) \sim \hat{\rho}_{\lambda}} ( (B,f) \notin\mathcal{A}_n)) (1-\mathbf{1}_{\mathcal{B}^c_n})\Bigr]
 \geq 
 \mathbb{E} \Bigl[ 1-\mathbb{P}_{ (B,f) \sim \hat{\rho}_{\lambda}} ( (B,f) \notin\mathcal{A}_n) - \mathbf{1}_{\mathcal{B}^c_n}
	\Bigr]
\geq 
1-\varepsilon_n
	\end{align*}
	will lead to
	\begin{equation*}
	\mathbb{E} \Bigl[ \mathbb{P}_{ (B,f) \sim \hat{\rho}_{\lambda}} ( (B,f)\in\mathcal{E}_n) \Bigr] \geq 1-\varepsilon_n.
	\end{equation*}
	In order to do so, assume that we are on the set $\mathcal{B}_n$, and let $ (B,f) \in\mathcal{A}_n$. Then,
	\begin{align*}
	\alpha    \Bigl( R(B,f) - R^\star \Bigr)
	 \leq 
	\inf_{\rho} \left\{ \lambda\Bigl( \int r_n(B,f) \rho({\rm d}(B,f)) - r_n^\star \Bigr)  +    \mathcal{K}(\rho,\pi)
	+ \log\frac{2}{\varepsilon_n} \right\}
 \leq
	\inf_{\rho} \left\{ \beta \Bigl( \int R(B,f) \rho({\rm d}(B,f)) - R^\star \Bigr)  +   2 \mathcal{K}(\rho,\pi)
	+ 2 \log\frac{2}{\varepsilon_n} \right\}
	\end{align*}
	that is,
	$$
	R(B,f) - R^\star 
	\leq 
	\inf_{\rho } \frac{ \beta \left[\int Rd\rho -
		R^\star \right] + 2 \left[
		\mathcal{K}(\rho, \pi) + \log \frac{2}{\varepsilon} \right] } {
		\alpha  }
	$$
	We upper-bound the right-hand side exactly as in the proof of Theorem~\ref{thm2}, this gives
	$(B,f)\in\mathcal{E}_n$.
	
\end{proof}

\section*{Conflicts of interest/Competing interests}
The authors declare no potential conflict of interests.

\section*{Acknowledgments}
TTM is supported by the Norwegian Research Council, grant number 309960 through the Centre for Geophysical Forecasting at NTNU.

\end{document}